\documentclass[a4paper,12pt]{article}
\usepackage{amsmath}
\usepackage{cases}
\usepackage{mathrsfs}
\usepackage{amsfonts}
\usepackage{latexsym}
\usepackage{indentfirst}
\usepackage[dvips]{graphicx}

\makeatletter

\newcommand{\Rmnum}[1]{\expandafter\@slowromancap\romannumeral #1@}
\title{On the Main Signless Laplacian Eigenvalues of a Graph
\thanks{This work was supported by Hunan Provincial Natural Science Foundation of China
(09JJ6009) and the Program for Science and Technology Innovative
Research Team in Higher Educational Institution of Hunan Province.}}
\author{{\sc Hanyuan Deng}\thanks{Corresponding author:
hydeng@hunnu.edu.cn.},  He Huang\\
{\small College of Mathematics and Computer Science,}\\
{\small Hunan Normal University, Changsha, Hunan 410081, P. R.
China}}
\date{2011-7-1}
\begin{document}
\maketitle

\begin{abstract}
A signless Laplacian eigenvalue of a graph $G$ is called a main
signless Laplacian eigenvalue if it has an eigenvector the sum of
whose entries is not equal to zero. In this paper, we first give the
necessary and sufficient conditions for a graph with one main
signless Laplacian eigenvalue or two main signless Laplacian
eigenvalues, and then characterize the trees and unicyclic graphs
with exactly two main signless Laplacian eigenvalues, respectively.

{\bf Keywords}: Signless Laplacian eigenvalue; main eigenvalue;
tree; unicyclic graph.
\end{abstract}

\baselineskip=0.30in

\section{Introduction}

Let $M$ be a square matrix of order $n$, an eigenvalue $\lambda$ of
$M$ is said to be a main eigenvalue if the eigenspace
$\varepsilon(\lambda)$ of $\lambda$ is not orthogonal to the all-1
vector $\mathbf{j}$, i.e., it has an eigenvector the sum of whose
entries is not equal to zero. An eigenvector $\mathbf{x}$ is a main
eigenvector if $\mathbf{x^{T}j}\neq 0$. Specially, if $M=A$ is the
$(0,1)$-adjacency matrix of a graph $G$, then the main eigenvalues
of $A$ are said to be main eigenvalues of $G$. A graph with exactly
one main eigenvalue is regular. Cvetkovi\'c~\cite{cv} proposed the
problem of characterizing graphs with exactly k main eigenvalues,
$k>1$. Hagos~\cite{ha} gave a characterization of graphs with
exactly two main eigenvalues. Recently, Hou and Zhou~\cite{hz}
characterized the tree with exactly two main eigenvalues. Hou and
Tian~\cite{ht} determined all connected unicyclic graphs with
exactly two main eigenvalues. Zhu and Hu~\cite{zhu} characterized
all connected bicyclic graphs with exactly two main eigenvalues.
Rowlinson~\cite{pr} surveyed results relating main eigenvalues and
main angles to the structure of a graph, and discussed graphs with
just two main eigenvalues in the context of measures of irregularity
and in the context of harmonic graphs.

In this paper, we assume that $G$ is a simple connected graph, and
will consider the main eigenvalues of the signless Laplacian matrix
$Q$ of $G$, where $Q=D+A$ and $D$ is the diagonal matrix of vertex
degrees. The main eigenvalues of $Q$ is said to be the main signless
Laplacian eigenvalues of $G$. The signless Laplacian appears very
rarely in published papers before 2003. Recently the signless
Laplacian has attracted the attention of researchers, see, e.g.
~\cite{dc1}~\cite{dc2}~\cite{dc3}~\cite{dc4}~\cite{dc5}. Here, we
will give the necessary and sufficient conditions for a graph with
one main signless Laplacian eigenvalue or two main signless
Laplacian eigenvalues, and then characterize the trees and unicyclic
graphs with exactly two main signless Laplacian eigenvalues,
respectively.

\section{The graphs with one or two main signless Laplacian eigenvalues}

In this section, we will show that a graph with exactly one main
signless Laplacian eigenvalue is regular, and give a
characterization of graphs with exactly two main signless Laplacian
eigenvalues.

Note that if $G$ is a simple connected graph with signless Laplacian
matrix $Q$, then there is an eigenvector $\mathbf{x}>0$ of the
largest eigenvalue $\mu_1$ of $Q$ such that
$Q\mathbf{x}=\mu_1\mathbf{x}$, and $x^T\mathbf{j}\neq0$ by the
Perron-Frobenius theorem. This shows that the largest eigenvalue
$\mu_1$ of $Q$ is a main signless Laplacian eigenvalue. So, $G$ has
at least one main signless Laplacian eigenvalue.

The following result gives a characterization of graphs with exactly
one main signless Laplacian eigenvalue.

{\bf Theorem 1}. A graph $G$ with exactly one main signless
Laplacian eigenvalue if and only if $G$ is regular.

{\bf Proof}. If $G$ is $k$-regular, then $Q\mathbf{j}=2k\mathbf{j}$.
This shows that $\mu_1=2k$ is an eigenvalue of $Q$ with an
eigenvector $\mathbf{j}$. Since $Q$ is a non-negative irreducible
symmetric matrix, $\mu_1=2k$ is the largest eigenvalue of $Q$ with
the multiplicity $1$ by the Perron-Frobenius theorem. And the
eigenvectors of other eigenvalues of $Q$ are orthogonal with
$\mathbf{j}$. So, $Q$ has exactly one main eigenvalue.

If $G$ has exactly one main signless Laplacian eigenvalue, then the
largest eigenvalue $\mu_1$ is the unique main eigenvalue of $Q$. Let
$\xi$ be a eigenvector of $\mu_1$, $\mathbf{V_1}=\varepsilon(\mu_1)$
the eigenspace of $\mu_1$, then $\mathbf{V_1}$ is the space spanning
by $\xi$. If $\mathbf{V_2}$ is the space spanning by eigenvectors of
all eigenvalues of $Q$ different from $\mu_1$, and $\mathbf{V_3}$ is
the space spanning by $\mathbf{j}$, then $dim(\mathbf{V_2})=n-1$ and
$dim(\mathbf{V_3})=1$. Since $Q$ is a real symmetric matrix,
$\mathbf{V_1}$ is the orthogonal complement of $\mathbf{V_2}$. And
$\mathbf{V_3}$ is also the orthogonal complement of $\mathbf{V_2}$
since $\mu_1$ is the unique main eigenvalue of $Q$. So,
$\mathbf{V_1}=\mathbf{V_3}$, and $\xi=a\mathbf{j}$ for some real
$a\neq 0$. From $Q\xi=\mu_1\xi$, the row sums of $Q$ are equal, and
$G$ is regular.  \hfill $\Box$

Now, we discuss the characterization of graphs with exactly two main
signless Laplacian eigenvalues.

For any positive semi-definite matrix $M$ of order $n$, all its
eigenvalues are non-negative. Let $\mu_1>\mu_2>\cdots>\mu_r$ be the
eigenvalues of $M$ with multiplicities $n_1,n_2,\cdots,n_r$,
respectively, where $n_1+n_2+\cdots+n_r=n$.
$\{\xi_{i1},\xi_{i2},\cdots,\xi_{in_i}\}$ is a standard and
orthogonal basis of the eigenspace $\varepsilon(\mu_i)$,
$i=1,2,\cdots,r$. And
$P=[\xi_{11},\cdots,\xi_{1n_1},\xi_{21},\cdots,\xi_{2n_2},\cdots,\xi_{r1},\cdots,\xi_{rn_r}]$,
$P_i=[0,\ldots,\xi_{i1},\ldots,\xi_{in_i},0,\ldots,0]$,
$\mathbf{j}=[1,1,\ldots,1]^T$, then
\begin{displaymath}
P^TMP= \left( \begin{array}{cccc}
\mu_1I_{n_1\times n_1} & 0 & \ldots & 0\\
0 & \mu_2I_{n_2\times n_2} & \ldots & 0\\
\vdots & \vdots  &  \ddots & \vdots\\
0 & 0 & \ldots & \mu_rI_{n_r\times n_r}
\end{array}\right)
\end{displaymath}

Let
\begin{displaymath}
E_i= \left( \begin{array}{ccccc}
0 & \ldots & \ldots & \ldots & 0\\
\vdots & \ddots & \vdots & \vdots & \vdots\\
\vdots &\vdots & I_{n_i\times n_i} & \vdots & \vdots\\
\vdots &\vdots &\vdots & \ddots &\vdots\\
 0 & \ldots & \ldots &\ldots & 0
\end{array}\right)
\end{displaymath}

$Q_i=PE_iP^T$, then
\begin{equation}\label{eq:Q}
Q_i\mathbf{j}=PE_iP^T\mathbf{j}=P_iP_i^T\mathbf{j}
\end{equation}

and $M$ has the spectral decomposition
$$M=\mu_1Q_1+\mu_2Q_2+\cdots+\mu_rQ_r$$
where
\begin{displaymath}
Q_iQ_j=\left\{\begin{array}{ll}
\mathbf{0}    & \textrm{$i\neq j$}\\
Q_i  & \textrm{$i=j$}
\end{array}\right.
\end{displaymath}
And, for any polynomial
$f(x)=a_0x^n+a_1x^{n-1}+\cdots+a_{n-1}x+a_n$,
\begin{align}
\notag f(M)&=a_0(M)^n+a_1(M)^{n-1}+\ldots+a_{n-1}M+a_nI\\
\notag       &=a_0\sum_{i=1}^r\mu_i^nQ_i+a_1\sum_{i=1}^r\mu_i^{n-1}Q_i+\cdots+a_{n-1}\sum_{i=1}^r\mu_iQ_i+a_n\sum_{i=1}^rQ_i\\
             &=\sum_{i=1}^rf(\mu_i)Q_i
\end{align}

{\bf Lemma 2}. Let $\mu_1$, $\mu_2$, $\ldots$, $\mu_t$ ($1\leq t\leq
r$) be the main eigenvalues of a positive semi-definite matrix $M$
of order $n$, and $m(x)=(x-\mu_1)(x-\mu_2)\cdots(x-\mu_t)$, then

(i) $m(M)\mathbf{j}=0$;

(ii) If $f(x)$ is a polynomial with real coefficients and
$f(M)\mathbf{j}=0$, then $m(x)|f(x)$.

{\bf Proof}. (i) From (2), we know that
$m(M)=\sum\limits_{i=1}^rm(\mu_i)Q_i=\sum\limits_{i=t+1}^rm(\mu_i)Q_i$.
And $$m(M)\mathbf{j}=\sum_{i=t+1}^rm(\mu_i)Q_i\mathbf{j}$$

Since $\mu_{t+1}, \cdots, \mu_r$ are not the main eigenvalues of
$M$, and from (\ref{eq:Q}), $Q_i\mathbf{j}=0$ for $i=t+1, \cdots,
r$. So,
$m(M)\mathbf{j}=\sum\limits_{i=t+1}^rm(\mu_i)Q_i\mathbf{j}=0$.

(ii) From (2) and $Q_i\mathbf{j}=0$ for $i=t+1, \cdots, r$,
$$f(M)\mathbf{j}=\sum\limits_{i=1}^rf(\mu_i)Q_i\mathbf{j}=\sum\limits_{i=1}^tf(\mu_i)Q_i\mathbf{j}.$$

Since $f(M)\mathbf{j}=0$,
$\sum\limits_{i=1}^tf(\mu_i)Q_i\mathbf{j}=0$. For $k=1, 2, \ldots,
t$, we have
$Q_k(\sum\limits_{i=1}^tf(\mu_i)Q_i\mathbf{j})=f(\mu_k)Q_k\mathbf{j}=0$.
So, $f(\mu_k)=0$ for $k=1, 2, \ldots, t$ and $m(x)|f(x)$. \hfill
$\Box$

A number $\alpha$ is an algebraic integer if there is a monic
polynomial $f(x)$ with integral coefficients such that
$f(\alpha)=0$.

{\bf Lemma 3}(~\cite{D.S}). $\alpha\in\mathbb{Q}$ is an algebraic
integer if and only if $\alpha$ is an integer.

{\bf Lemma 4}(~\cite{D.S}). If $\alpha$ and $\beta$ are algebraic
integers, then $\alpha\pm\beta$ and $\alpha\beta$ are also algebraic
integers.

{\bf Theorem 5}. Let $G$ be non-regular, then $G$ has exactly two
main signless Laplacian eigenvalues $\mu_1$ and $\mu_2$ if and only
if $(Q-\mu_1I)(Q-\mu_2I)\mathbf{j}=0$.

{\bf Proof}. Let $\mu_1, \ldots, \mu_t$ be the main eigenvalues of
$Q$, and $m(x)=(x-\mu_1)\cdots(x-\mu_t)$.

If $(Q-\mu_1I)(Q-\mu_2I)\mathbf{j}=0$, then $f(Q)\mathbf{j}=0$ for
$f(x)=(x-\mu_1)(x-\mu_2)$, and $m(x)|f(x)$ by Lemma 2. So,
$m(x)=(x-\mu_1)(x-\mu_2)$ or $(x-\mu_1)$ or $(x-\mu_2)$, and $t\leq
2$. But $G$ is non-regular, $t=2$ from Theorem 1.

If $G$ has exactly two main signless Laplacian eigenvalues $\mu_1$
and $\mu_2$, then $(Q-\mu_1I)(Q-\mu_2I)\mathbf{j}=0$ from Lemma 2.
\hfill $\Box$

In the following, we give an alternative characterization of graphs
with exactly two main signless Laplacian eigenvalues.

In order to find all graphs with exactly two main eigenvalues, Hou
and Tian~\cite{ht} introduced a 2-walk $(a,b)$-linear graph. For a
graph $G$, the degree of vertex $v$ is denoted by $d(v)$, the number
of walks of length $2$ of $G$ starting at $v$ is
$s(v)=\sum\limits_{u\in N_G(v)}d(u)$, i.e., the sum of the degrees
of the vertices adjacent to $v$, where $N_G(v)$ is the set of all
neighbors of $v$ in $G$. A graph $G$ is called 2-walk $(a,
b)$-linear if there exist unique integer numbers $a$, $b$ with
$a^2-4b>0$ such that $s(v)=ad(v)+b$ holds for every vertex $v\in
V(G)$. Hagos~\cite{ha} showed that a graph $G$ has exactly two main
eigenvalues if and only if $G$ is 2-walk linear.

Like a 2-walk $(a, b)$-linear graph, we definite a 2-walk $( a,
b)$-parabolic graph. A graph $G$ is called 2-walk $(a, b)$-parabolic
if there are uniquely a positive integer $a$ and a non-negative
integer $b$ with $a^2-8b>0$ such that $s(v)=-d^2(v)+ad(v)-b$ holds
for every vertex $v\in V(G)$.

{\bf Theorem 6}. A graph $G$ has exactly two main signless Laplacian
eigenvalues if and only if $G$ is a 2-walk $(a, b)$-parabolic graph.

{\bf Proof}. If $G$ is a 2-walk $(a, b)$-parabolic graph, then there
are uniquely a positive integer $a$ and a non-negative integer $b$
such that $a^2-8b>0$ and $s(v)=-d^2(v)+ad(v)-b$ for any $v\in
V(G)=\{v_1,v_2,\cdots,v_n\}$. So, $s(v_i)+d^2(v_i)-ad(v_i)+b=0$, and
\begin{align}
\quad \quad \notag \frac{1}{2}(A+D)^2\mathbf{j}-aA\mathbf{j}+b\mathbf{j}&=\mathbf{0} \\
\notag
\frac{1}{2}{Q}^2\mathbf{j}-\frac{1}{2}aL^{+}\mathbf{j}+b\mathbf{j}&=\mathbf{0}\\
\notag {Q}^2\mathbf{j}-aQ\mathbf{j}+2b\mathbf{j}&=\mathbf{0}
\end{align}

Let $f(x)=x^2-ax+2b$, then $f(Q)\mathbf{j}=\mathbf{0}$, and $f(x)=0$
has two real roots since $a^2-8b>0$. And $G$ is non-regular since
one has $s(v)=-d^2(v)+2kd(v)-0$ and $s(v)=-d^2(v)+(2k+1)d(v)-k$ for
a $k$-regular graph, i.e., $(a, b)=(2k, 0)$ or $(2k+1, k)$ is not
unique. From Theorem 5, $G$ has exactly two main signless Laplacian
eigenvalues.

On the other hand, if $G$ has exactly two main signless Laplacian
eigenvalues $\mu_1$ and $\mu_2$, then by Theorem 5,
$$({Q}^2-(\mu_1+\mu_2)Q+\mu_1\mu_2I)\mathbf{j}=\mathbf{0}$$
i.e.,
$$(D+A)^2\mathbf{j}-(\mu_1+\mu_2)(D+A)\mathbf{j}+\mu_1\mu_2\mathbf{j}=\mathbf{0}.$$
So, $d^2(v)+s(v)-(\mu_1+\mu_2)d(v)+\frac{\mu_1\mu_2}{2}=0$ for all
$v\in V(G)$. Let $\mu_1+\mu_2=a$ and $\mu_1\mu_2=2b$, then
$s(v)=-d^2(v)+ad(v)-b$, and $a>0$, $b\geq 0$ and $a^2-8b>0$ since
$\mu_1\neq\mu_2$ are the eigenvalues of the positive semi-definite
matrix $Q$. Note that $G$ is non-regular by Theorem 1, there are $u,
v\in V(G)$ such that $d(u)\neq d(v)$. From $s(u)=-d^2(u)+ad(u)-b$
and $s(v)=-d^2(v)+ad(v)-b$, we have
\begin{displaymath}
a=\frac{s(u)-s(v)}{d(u)-d(v)}+d(u)+d(v)
\end{displaymath}
\begin{equation}\label{eq:a4}
b=\frac{s(u)-s(v)}{d(u)-d(v)}d(v)+d(u)d(v)-s(v)
\end{equation}
and $a, b$ are rational numbers and unique. Because $\mu_1, \mu_2$
are the roots of monic polynomial $det(\lambda I-Q)=0$ with integral
coefficients, $\mu_1, \mu_2$ are algebraic integers. By Lemmas 3 and
4, $a, b$ are integers. \hfill $\Box$

\section{Trees with exactly two main signless Laplacian eigenvalues}

In this section, we will determine all trees with exactly two main
signless Laplacian eigenvalues.

Let $G=(V,E)$ be a tree with $n\geq 3$ vertices and the maximum
degree $\Delta$. If $G$ has exactly two main signless Laplacian
eigenvalues, then from Theorem 6, there exist uniquely a positive
integer $a$ and a non-negative integer $b$ such that such that
$a^2-8b>0$ and
\begin{equation}\label{eq:s}
s(v)=-d^2(v)+ad(v)-b
\end{equation}
for any $v\in V(G)$.

{\bf Case 1}. $b=0$.

Let $v_1\in V$ with degree $d(v_1)=1$, and $v_2$ is its unique
adjacent vertex. Then $d(v_2)=s(v_1)=-1+a$ by (\ref{eq:s}), and
\begin{equation}\label{eq:a}
a=d(v_2)+1\leq\Delta+1
\end{equation}
Let $v_0\in V$ with degree $d(v_0)=\Delta$, then $\Delta=d(v_0)\leq
s(v_0)=-d^2(v_0)+ad(v_0)=-\Delta^2+a\Delta$, and $\Delta\leq a-1$,
i.e.,
\begin{equation}\label{eq:a1}
a\geq \Delta+1
\end{equation}
with equality if and only if $G$ is a star with the center $v_0$.
From (\ref{eq:a}) and (\ref{eq:a1}), we have $a=\Delta+1$. So,
$G=S_n$ is a star.

{\bf Case 2}. $b=1$.

Let $P_k=v_1v_2\ldots v_k$ is a longest path of $G$. Then
$d(v_2)=s(v_1)=-1+a-1=a-2$ by (\ref{eq:s}), and
\begin{align}
\notag \quad s(v_2)&=-d^2(v_2)+ad(v_2)-1 \\
\notag       &=-(a-2)^2+a(a-2)-1\\
\notag       &=2a-5
\end{align}
Since $P_k=v_1v_2\ldots v_k$ is a longest path of $G$, the adjacent
vertices of $v_2$ are pendant vertices except $v_3$. So,
$s(v_2)=\sum\limits_{v\in N_G(v_2)}d(v)=d(v_3)+d(v_2)-1$ and
\begin{align}
\notag d(v_3)&=s(v_2)-d(v_2)+1\\
\notag       &=(2a-5)-(a-2)+1\\
\notag       &=a-2=d(v_2)
\end{align}
By (\ref{eq:s}), we have $s(v_3)=s(v_2)=2a-5=d(v_2)+d(v_3)-1$. This
shows that the adjacent vertices of $v_3$ are pendant vertices
except $v_2$. So, $G=S_{\frac{n}{2},\frac{n}{2}}$ is a double star.

{\bf Case 3}. $b\geq 2$.

Let $P_k=v_1v_2\ldots v_k$ is a longest path of $G$, then
$d(v_2)\geq 2$. By (\ref{eq:s}), $d(v_2)=s(v_1)=-1+a-b$, and
$a-b\geq 3$.

Since $P_k=v_1v_2\ldots v_k$ is a longest path of $G$, the adjacent
vertices of $v_2$ are pendant vertices except $v_3$.
\begin{align}
\notag d(v_3)&=s(v_2)-(d(v_2)-1) \\
\notag       &=(-(a-b-1)^2+a(a-b-1)-b)-(d(v_2)-1)\\
\notag       &=(-(a-b-1)^2+a(a-b-1)-b)-(a-b-2)\\
\notag       &=ab-b^2-2b+1
\end{align}
\begin{align}
\mbox{and} \qquad \notag d(v_3)-d(v_2)&=(ab-b^2-2b+1)-(a-b-1)\\
\notag              &=ab-b^2-b-a+2\\
\notag              &=(b-1)(a-b-2)>0
\end{align}
So, $d(v_3)>d(v_2)=a-b-1\geq 2$. And no pendant vertex is adjacent
to $v_3$; Otherwise, let $u$ be a pendant vertex adjacent to $v_3$.
Then $d(v_3)=s(u)=-1+a-b$ by (\ref{eq:s}), contradicting with
$d(v_3)>a-b-1$.

For any $x\in N_{G}(v_3)\backslash \{v_2,v_4\}$, since $x$ is not a
pendant vertex, there is $y\in V(G)\backslash\{v_3\}$ such that
$xy\in E(G)$, and $y$ is a pendant vertex by the longest path
$P_k=v_1v_2\ldots v_k$. Then
\begin{displaymath}
d(x)=s(y)=-1+a-b=d(v_2),  \,\,\,\,\forall x\in N_{G}(v_3)\backslash
\{v_2,v_4\}
\end{displaymath}
So,
\begin{align}
\quad \notag  s(v_3)=\sum_{z\in N_G(v_3)}d(z)&=\sum_{x\in
N_G(v_3)\backslash\{v_2,v_4\}}d(x)+d(v_2)+d(v_4)\\
\notag       &=(d(v_3)-2)d(v_2)+d(v_2)+d(v_4)\\
\mbox{and}\qquad \notag d(v_4)&=s(v_3)-(d(v_3)-1)d(v_2)
\end{align}

By (\ref{eq:s}), $s(v_3)=-d^2(v_3)+ad(v_3)-b$. Note that
$d(v_2)=a-b-1$ and $d(v_3)=ab-b^2-2b+1$,
\begin{align}
\notag d(v_4)&=-d^2(v_3)+ad(v_3)-b-(d(v_3)-1)d(v_2)\\
\notag       &=d(v_3)(-d(v_3)+a-d(v_2))-b+d(v_2)\\
\notag       &=d(v_3)(-ab+b^2+3b)+a-2b-1\\
\notag       &=d(v_3)b(b-a+3)+(a-2b-1)
\end{align}

If $a-b=3$, then $d(v_4)=a-2b-1=2-b\leq0$; If $a-b\geq 4$, then
$d(v_4)=d(v_3)b(b-a+4)-bd(v_3)+(a-2b-1)\leq
-bd(v_3)+(a-2b-1)=-b(ab-b^2-2b+1)+(a-2b-1)=-a(b^2-1)+b^3+2b^2-3b-1\leq
-(b+4)(b^2-1)+b^3+2b^2-3b-1=-2b^2-2b+3<0$. This is impossible.

On the other hand, it is easy to check that $G=S_n$ and
$G=S_{\frac{n}{2},\frac{n}{2}}$ are 2-walk $(n,0)$-parabolic graph
and $(\frac{n}{2}+2,1)$-parabolic graph, respectively.

From above, we have

{\bf Theorem 7}. A tree with $n\geq 3$ vertices has exactly two main
signless Laplacian eigenvalues if and only if $G$ is the star $S_n$
or the double star $S_{\frac{n}{2},\frac{n}{2}}$.

It was showed in~\cite{hz} that the trees with $n\geq 3$ vertices
has exactly two main eigenvalues (of adjacent matrix) are $S_n$,
$S_{\frac{n}{2},\frac{n}{2}}$ and $T_a$. But from Theorem 7, we know
that $T_a$ is not a tree with exactly two main signless Laplacian
eigenvalues, where $T_a$ ($a\geq 2$) is defined in~\cite{hz} to be
the tree with one vertex $v$ of degree $a^2-a+1$ while every
neighbor of $v$ has degree $a$ and all remaining vertices are
pendant.

\section{Unicyclic graphs with exactly two main signless Laplacian eigenvalues}

In this section, we will determine all unicyclic graphs with exactly
two main signless Laplacian eigenvalues.

The unique unicyclic graph with $n$ vertices and the minimum degree
$\delta\geq 2$ is the cycle $C_n$, and it is regular. By Theorem 1,
it has exactly one main signless Laplacian eigenvalues. So, we only
need to consider the unicyclic graphs with the minimum degree
$\delta=1$.

{\bf Remark 1}. If $G$ is a 2-walk $(a, b)$-parabolic graph with
$\delta(G)=1$, then $a-b\geq 3$ since there is a pendent vertex $x$
with the only incident edge $xy$ in $G$ and $d(y)=s(x)=-1+a-b\geq
2$.

Let $\mathscr{G}_{a,b}=\{G: G \mbox{ is a 2-walk } (a,
b)-\mbox{parabolic unicyclic graph with}$ $\delta(G)=1\}$, and for
each $G\in \mathcal{G}_{a,b}$, let $G_0$ be the graph obtained from
$G$ by deleting all pendant vertices. If $v\in V(G_0)$, we use
$d_{G_0}(v)$ to denote the degree of the vertex $v$ in $G_0$.

{\bf Lemma 8}. If $G\in\mathscr{G}_{a,b}$ and $v\in V(G_0)$, then
$d(v)=d_{G_0}(v)$ or $d(v)=a-b-1$.

{\bf Proof}. If there is a pendant $x$ adjacent to $v$ in $G$, then
$d(v)=s(x)=a-b-1$ by (\ref{eq:s}). Otherwise, $d(v)=d_{G_0}(v)$.
\hfill $\Box$

{\bf Lemma 9}. If $G\in\mathscr{G}_{a, b}$, then (i)
$\delta(G_0)\geq 2$; (ii) $a-b\geq 4$ and $a\geq 5$.

{\bf Proof}. (i) If $\delta(G_0)=1$, then there is $y\in V(G_0)$
such that $d_{G_0}(y)=1$, and there must exist a pendant vertex $x$
adjacent to $y$ in $G$. By (\ref{eq:s}), $d(y)=s(x)=-1+a-b$, this
shows that there are $a-b-2$ pendant vertices adjacent to $y$ in
$G$. Let $z$ be the unique non-pendant vertex adjacent to $y$ in
$G$, then $s(y)=\sum\limits_{w\in N_G(y)}d(w)=d(z)+(a-b-2)$. By
(\ref{eq:s}), we know
\begin{align}
\notag s(y)&=-d^2(y)+ad(y)-b\\
\notag     &=-(a-b-1)^2+a(a-b-1)-b
\end{align}
and $d(z)=s(y)-(a-b-2)=-(a-b-1)^2+a(a-b-1)-b-(a-b-2)$,
\begin{equation}\label{eq:8}
\mbox{ i.e., } \qquad \quad \qquad \quad d(z)=ab-b^2-2b+1
\end{equation}
So,
\begin{equation}\label{eq:9}
d(z)-d(y)=(ab-b^2-2b+1)-(a-b-1)=(b-1)(a-b-2)
\end{equation}

(I) If $b=0$, then by (\ref{eq:8}), $d(z)=1$. This is impossible
since $z$ is a non-pendant vertex adjacent to $y$ in $G$.

(II) If $b=1$, then by (\ref{eq:8}) and (\ref{eq:9}), $d(z)=a-2$ and
$d(y)=d(z)=a-2$.

From (\ref{eq:s}), $s(z)=-(a-2)^2+a(a-2)-1=2a-5=d(y)+d(z)-1$. This
shows that all the vertices adjacent to $z$, except $y$, are pendant
vertices. So, $G$ is a double star with the centers $z$ and $y$.
This is impossible since $G\in\mathscr{G}_{a, b}$.

(III) If $b\geq 2$, then no pendant vertex is adjacent to $z$ in
$G$; Otherwise, $d(z)=s(u)=-1+a-b\geq 2$, where $u$ is a pendant
vertex is adjacent to $z$. This implies that $d(z)=d(y)$ and
$a-b\geq 3$, contradicting with (\ref{eq:9}).

By (\ref{eq:s}), we have $s(z)=-d^2(z)+ad(z)-b$. And
\begin{align}
\notag s(z)&=\sum\limits_{w\in N_{G}(z)}d(w)=\sum_{w\in N_{G_0}(z)}d(w)\\
\notag     &=d(y)+\sum_{w\in N_{G_0}(z)\backslash \{y\}}d(w)\\
\notag     &\geq d(y)+2(d(z)-1)\\
\notag     &=a-b-1+2(d(z)-1)
\end{align}
So, $-d^2(z)+ad(z)-b\geq a-b-1+2(d(z)-1)$
\begin{displaymath}
d^2(z)-(a-2)d(z)+a-3\leq 0
\end{displaymath}
and
\begin{displaymath}
1\leq d(z)\leq a-3
\end{displaymath}

By (\ref{eq:8}), $ab-b^2-2b+1\leq a-3$, i.e., $ b^2+(2-a)b+a-4\geq
0$. Then
$$b\leq\frac{(a-2)-\sqrt{(a-4)^2+4}}{2}\,\,\,\, \mbox{ or }\,\,\,
b\geq\frac{(a-2)+\sqrt{(a-4)^2+4}}{2}$$

From Theorem 6, $a^2>8b\geq 16$, i.e., $a>4$, and
$\frac{(a-2)-\sqrt{(a-4)^2+4}}{2}<\frac{(a-2)-(a-4)}{2}=1$. So,
$b\geq\frac{(a-2)+\sqrt{(a-4)^2+4}}{2}$. But
$$2\leq d(y)=
a-b-1\leq
a-\frac{(a-2)+\sqrt{(a-4)^2+4}}{2}-1=\frac{a-\sqrt{(a-4)^2+4}}{2}.$$
We have $a-4\geq\sqrt{(a-4)^2+4}$. This is impossible.

Summarizing (I)-(III) above, we have $\delta(G_0)\geq 2$.

(ii) Because $G\in\mathscr{G}_{a, b}$, $\delta(G)=1$. There is a
pendent vertex $x$ and the only edge $xy$ incident with $x$ in $G$.
$d(y)=s(x)=a-b-1$ by (\ref{eq:s}). From (i), $d(y)\geq
d_{G_0}(y)+1\geq\delta(G_0)+1\geq3$. So, $a-b\geq 4$.

Since \begin{align}
\qquad \quad \notag s(y)&=d(y)-d_{G_0}(y)+\sum_{w\in N_{G_0}(y)}d(w)\\
\notag     &\geq d(y)-d_{G_0}(y)+2d_{G_0}(y)\\
\notag     &\geq a-b-1+2=a-b+1,
\end{align}
by (\ref{eq:a4}), we have
\begin{align}
\quad  \notag a&=\frac{s(y)-s(x)}{d(y)-d(x)}+d(y)+d(x)\\
\notag  &\geq \frac{(a-b+1)-(a-b-1)}{a-b-2}+(a-b-1)+1\\
\notag  &=\frac{2}{a-b-2}+a-b>a-b\geq 4.
\end{align}
So, $a\geq 5$ since $a$ is a integer from Theorem 6.  \hfill $\Box$

In the following, we determine all unicyclic graphs with exactly two
main signless Laplacian eigenvalues.

Let $G_1$ be the unicyclic graph with $n$ vertices obtained by
attaching $k\geq 1$ pendant vertices to each vertex of a cycle with
length $r$, where $n=(k+1)r$. It was showed in~\cite{ht} that $G_1$
is the only connected graph with exactly two main eigenvalue (of
adjacent matrix). $G_2$ is the unicyclic graph with $n$ vertices
obtained from the cycle $u_1u_2\cdots u_{3t}$ by attaching one
pendant vertices to the vertex $u_{3s+1}$ for $s=0,1,\cdots,t-1$,
where $n=4t$ (see Figure 1).

{\bf Theorem 10}. Let $G$ be a unicyclic graphs with $n$ vertices
different from the cycle $C_n$. $G$ has exactly two main signless
Laplacian eigenvalues if and only if $G$ is isomorphic to one of
graphs $G_1$ and $G_2$, both given in Figure 1.

{\bf Proof}. First, it is easy to check that $G_1$ is a 2-walk
$(k+5, 2)$-parabolic graph and $G_2$ is a 2-walk $(5,1)$-parabolic
graph.

Next, because $G\in \mathscr{G}_{a,b}$ is unicyclic, $G_0$ is a
cycle from Lemma 9. Let $G_0=C_r=u_1u_2\ldots u_ru_1$, then
$d(u_i)\in\{d_{G_0}(u_i), a-b-1\}=\{2, a-b-1\}$ from Lemma 8, where
$1\leq i\leq r$.

(i) If $d(u_1)=d(u_2)=\cdots=d(u_r)=a-b-1$, then $u_i$ has $k=a-b-3$
pendant vertices for $1\leq i\leq r$. So, $G\cong G_1$.

(ii) If there is $u_i\in V(G_0)$ such that $d(u_i)=2$ for some
$i\in\{1,2,\cdots,r\}$, then by (\ref{eq:s}), we have
\begin{equation}\label{eq:10}
d(u_{i-1})+d(u_{i+1})=s(u_i)=-d^2(u_i)+ad(u_i)-b=-4+2a-b
\end{equation}
Without loss of generality, we assume that $d(u_{i+1})\geq
d(u_{i-1})$.  From Lemma 9(ii),
$d(u_{i-1})+d(u_{i+1})=s(u_i)=(a-b-4)+a\geq 5$. Since $d(u_{i-1}),
d(u_{i+1})\in\{d_{G_0}(u_i), a-b-1\}=\{2, a-b-1\}$, we have
$d(u_{i+1})=d(u_{i-1})=a-b-1$, or $d(u_{i+1})=a-b-1$ and
$d(u_{i-1})=2$.

If $d(u_{i+1})=d(u_{i-1})=a-b-1$, then $b=2$ by (\ref{eq:10}). From
(\ref{eq:s}), $s(u_{i-1})=-(a-b-1)^2+a(a-b-1)-b=3a-11$. And
$$s(u_{i-1})=d(u_{i-2})+d(u_i)+(d(u_{i-1})-d_{G_0}(u_{i-1}))$$
i.e.,
$d(u_{i-2})=s(u_{i-1})-d(u_i)-d(u_{i-1})+d_{G_0}(u_{i-1})=(3a-11)-2-(a-b-1)+2=2a-8$.
From Lemma 8,
$d(u_{i-2})=2a-8\in\{a-b-1,d_{G_0}(u_{i-2})\}=\{a-3,2\}$. We have
$a=5$, and $a-b=3<4$, contradicting with Lemma 9(ii).

If $d(u_{i+1})=a-b-1$ and $d(u_{i-1})=2$, then $a=5$ from
(\ref{eq:10}). By (\ref{eq:s}),
$s(u_{i+1})=-(a-b-1)^2+a(a-b-1)-b=4+2b-b^2$. And
$$s(u_{i+1})=d(u_{i+2})+d(u_i)+(d(u_{i+1})-2)$$
i.e.,
$d(u_{i+2})=s(u_{i+1})-d(u_i)-d(u_{i+1})+2=(4+2b-b^2)-2-(a-b-1)+2=3b-b^2$.
From Lemma 8,
$d(u_{i+2})=3b-b^2\in\{a-b-1,d_{G_0}(u_{i-2})\}=\{4-b,2\}$. We have
$b=1$ or $b=2$. And $a-b\geq 4$ from Lemma 9. So, $b=1$ and
$d(u_{i+1})=3$, $d(u_{i+2})=2$.

By (\ref{eq:s}) again,

$d(u_{i+3})+d(u_{i+1})=s(u_{i+2})=-d^2(u_{i+2})+5d(u_{i+2})-1=5$,
and $d(u_{i+3})=5-d(u_{i+1})=2$;

$d(u_{i+4})+d(u_{i+2})=s(u_{i+3})=-d^2(u_{i+3})+5d(u_{i+3})-1=5$,
and $d(u_{i+4})=5-d(u_{i+2})=3$;

$d(u_{i+5})+d(u_{i+3})+1=s(u_{i+4})=-d^2(u_{i+4})+5d(u_{i+4})-1=5$,
and $d(u_{i+5})=5-d(u_{i+3})-1=2$.

Continuing like this, we have
\begin{displaymath}
d(u_k) = \left\{ \begin{array}{ll}
2, & \mbox{ $k-i\equiv 0,2(mod3)$}; \\
3, & \mbox{ $k-i\equiv 1(mod3)$}
\end{array} \right.
\end{displaymath}
So, $r\equiv 0(mod3)$ and $G\cong G_2$. \hfill $\Box$

The results on main sigenless Laplacian eigenvalues presented in
Section 2 are useful to the problem of characterizing graphs with a
given number of main sigenless Laplacian eigenvalues. And Theorems 7
and 10 show that the set of graphs with a given number of main
sigenless Laplacian eigenvalues is not identical with the set of
graphs a given number of main eigenvalues (of adjacent matrix).

\end{document}